\documentclass[review]{elsarticle}

\usepackage{lineno,hyperref}
\modulolinenumbers[5]
\usepackage[b5paper]{geometry}	    
\usepackage{amsfonts,amsmath,latexsym,amssymb} 
\usepackage{theorem}               
\usepackage{mathrsfs,upref}        
\usepackage{mathptmx}		 
	                            
\journal{Journal of \LaTeX\ Templates}

\newtheorem{theorem}{Theorem}

\theoremstyle{definition}

\newtheorem{example}{Example}

\begin{document}

\begin{frontmatter}

\title{Extensions of the Schur majorisation inequalities}

\author{Rajendra Bhatia}
\address{Ashoka University, Sonepat, Haryana, 131029, India}
\author{Rajesh Sharma}
\address{Department of Mathematics and Statistics, H.P. University, Shimla-5, India}

\begin{abstract}
Let $\lambda_{j}$ and $a_{jj}$, $1\leq j\leq n,$ be the eigenvalues and the diagonal entries of a Hermitian matrix $A$, both enumerated in the increasing order. We prove some inequalities that are stronger than the Schur majorisation inequalities $\sum_{j=1}^{r}\lambda_{j}\leq \sum_{j=1}^{r}a_{jj},$ $1\leq r\leq n.$
\end{abstract}


\end{frontmatter}
\section{Introduction}
Let $A$ be an $n\times n$ complex Hermitian matrix. Let the eigenvalues and the diagonal entries of $A$ both be enumerated in increasing order as
\begin{equation}
\label{aa1}\lambda_{1}(A)\leq \lambda_{2}(A)\leq\dots\leq\lambda_{n}(A),   \tag{1.1}
\end{equation}
and
\begin{equation}
\label{aa2}a_{11}\leq a_{22}\leq\dots\leq a_{nn}, \tag{1.2}
\end{equation}
respectively. We then have
\begin{equation}
\label{aa3}\lambda_{1}(A)\leq a_{11} \quad\text{and}\quad\lambda_{n}(A)\geq a_{nn}.  \tag{1.3}
\end{equation}
These inequalities are included in the Schur majorisation inequalities that say: for every $1\leq r\leq n$
\begin{equation}
\label{aa4} \sum_{j=1}^{r}\lambda_{j}(A)\leq \sum_{j=1}^{r}a_{jj}, \tag{1.4}
\end{equation}
with equality in the case $r=n.$ These inequalities are of fundamental importance in matrix analysis and have been the subject of intensive work.
See, e.g. Bhatia \cite{bha97}, Horn and Johnson \cite{hor13} and Marshal and Olkin \cite{mar79}.\\
In this note we obtain some inequalities that are
stronger than \eqref{aa3} and \eqref{aa4}. These give estimates of eigenvalues in terms of quantities easily computable from the entries of $A.$\par Given the $n\times n$ Hermitian matrix $A=\left[ a_{ij}\right],$ let 
\begin{equation}
\label{aa5} r_{i}=\sum_{j\neq i} \left\vert\ a_{ij}\right\vert,\quad 1\leq i\leq n \tag{1.5}
\end{equation}
and
\begin{equation}
\label{aa6} q_{i}= \sum_{j\neq i} \left\vert\ a_{ij}\right\vert^{2},\quad 1\leq i\leq n. \tag{1.6}
\end{equation}
A permutation similarity does not change either the eigenvalues or the diagonal entries of $A.$ Nor does it change the quantities $r_{i}$ and $q_{i}.$ We assume that such a permutation similarity has been performed and the ordering \eqref{aa2} for diagonal entries has been achieved. To rule out trivial cases, we assume that $A$ is not a diagonal matrix.\\Our first theorem is a strengthening of the inequalities \eqref{aa3}.
\begin{theorem}
\label{th1} For every $n\times n$ Hermitian matrix $A$, we have
\begin{equation}
\label{aa7} \lambda_{1}(A)\leq a_{11}-\frac{q_{1}}{\max_{i}(a_{ii}+r_{i})-a_{11}}, \tag{1.7}
\end{equation}

\begin{equation}
\label{aa8} \lambda_{n}(A)\geq a_{nn}+\frac{q_{n}}{a_{nn}-\min_{i}(a_{ii}-r_{i})}.  \tag{1.8}
\end{equation}
\end{theorem}
The next two theorems give inequalities stronger than (1.4).
\begin{theorem}
\label{th2}Let $A$ be an $n\times n$ Hermitian matrix. Then for $1\leq r\leq n-1$ and $r<t\leq n$, we have
\begin{equation}
\label{aa9} \sum_{i=1}^{r}\lambda_{i}(A)\leq \sum_{i=1}^{r}a_{ii}-\frac{\sum_{s=1}^{r} \left\vert\ a_{ts}\right\vert^{2}}{a_{tt}-\min_{i=1,...,r,t}\left(a_{ii}-\sum_{\substack{s=1\\s\neq i}}^{r+1} \left\vert\ a_{is}\right\vert\right)}.  \tag{1.9}
\end{equation}
\end{theorem}
\begin{theorem}
\label{th3}Let $A$ be an $n\times n$ Hermitian matrix. Then for $1\leq r\leq n-1$, $1\leq k\leq r$, and $r<t\leq n$, we have
\begin{equation}
\label{aa10} \sum_{i=1}^{r}\lambda_{i}(A)\leq \sum_{i=1}^{r}a_{ii}-\frac{\sqrt{(a_{tt}-a_{kk})^{2}+4|a_{tk}|^{2}}-(a_{tt}-a_{kk})}{2}.   \tag{1.10}
\end{equation}
\end{theorem}

\section{Proofs}
Our proofs rely upon two basic theorems of matrix analysis. Let $\mathbb{M}(n)$ be the algebra of all $n\times n$ complex matrices and let $\Phi :\mathbb{M}%
(n)\rightarrow \mathbb{M}(k)$ be a positive unital linear map,
\cite{bha07}. Then the Bhatia-Davis inequality \cite{bha00} says that for every Hermitian matrix $A$ whose spectrum is contained in the interval $\left[ m,M\right],$ we have
\begin{equation*}
 \label{bb1}\Phi \left( A^{2}\right)-\Phi \left( A\right)^{2}\leq\left(MI-\Phi \left( A\right)\right)\left(\Phi \left( A\right)-mI\right)
\leq\left(\frac{M-m}{2}\right)^{2}I.  \tag{2.1}
\end{equation*}
Cauchy's interlacing principle says that if $A_{r}$ is an $r\times r$ principal submatrix of $A,$ then
\begin{equation}
\label{bb2}\lambda_{j}(A)\leq \lambda_{j}(A_{r}),\quad 1\leq j\leq r.  \tag{2.2}
\end{equation}
See Chapter III of \cite{bha97} for this and other facts used here.
\subsection{Proof of Theorem 1}
Let $\varphi :\mathbb{M}%
(n)\rightarrow \mathbb{C}$ be a positive unital linear functional and let the eigenvalues of Hermitian element $A\in\mathbb{M}(n)$ be arranged as in \eqref{aa1}. From the first inequality \eqref{bb1}, we have
\begin{equation}
 \label{bb3}\varphi \left( A^{2}\right)-\varphi \left( A\right)^{2}\leq\left(\lambda_{n}(A)-\varphi \left( A\right)\right)\left(\varphi \left( A\right)-\lambda_{1}(A)\right). \tag{2.3}
\end{equation}
Suppose $\lambda_{n}(A)\neq \varphi(A)$. Then, from \eqref{bb3}, we have
\begin{equation}
\label{bb4} \lambda_{1}(A)\leq \varphi(A)- \frac{\varphi\left(A^{2}\right)-\varphi\left(A\right)^{2}}{\lambda_{n}(A)-\varphi\left(A\right)}. \tag{2.4}
\end{equation}
Further, by the Gersgorin disk theorem, we have
\begin{equation}
\label{bb5}\lambda_{n}(A)\leq \text{max}_{i}(a_{ii}+r_{i}).  \tag{2.5}
\end{equation}
Combining \eqref{bb4} and \eqref{bb5}, we get
\begin{equation}
\label{bb6}\lambda_{1}(A)\leq \varphi(A)-\frac{\varphi\left(A^{2}\right)-\varphi\left(A\right)^{2}}{\text{max}_{i}(a_{ii}+r_{i})-\varphi\left(A\right)}.  \tag{2.6}
\end{equation}
Choose $\varphi(A)=a_{11}$. Then, $\varphi$ is a positive unital linear functional and $\varphi \left( A^{2}\right)-\varphi \left( A\right)^{2}=q_{1}$. So, \eqref{bb6} yields \eqref{aa7}.\\Suppose $\lambda_{n}(A)=\varphi(A)=a_{11}.$ Then, from \eqref{aa2} and \eqref{aa3}, we have $a_{11}=a_{22}=\dots=a_{nn}$ and from \eqref{bb3}, $\varphi \left( A^{2}\right)-\varphi \left( A\right)^{2}=0.$ Therefore, $q_{i}=0$ for all $i=1,2,...,n$. But then $A$ is a scalar matrix.\\The inequality \eqref{aa8} follows on using similar arguments. The derivation requires lower bound of $\lambda_{n}(A)$ from \eqref{bb3} which is analogous to \eqref{bb4}, $\lambda_{1}(A)\geq\min_{i}\left(a_{ii}-r_{i}\right)$ and $\varphi \left( A\right)=a_{nn}$.   \qed
\subsection{Proof of Theorem 2}
The trace of $A$ is the sum of the eigenvalues of $A$. Therefore,
\begin{equation}
\label{bb7}\lambda_{n}(A)=\text{tr}A-\sum_{i=1}^{n-1}\lambda_{i}(A). \tag{2.7}
\end{equation}
Combining \eqref{aa8} and \eqref{bb7}, we find that
\begin{equation}
\label{bb8}\sum_{i=1}^{n-1}\lambda_{i}(A)\leq\sum_{i=1}^{n-1}a_{ii}-\frac{q_{n}}{a_{nn}-\min_{i}(a_{ii}-r)}.  \tag{2.8}
\end{equation}
Apply \eqref{bb8} to the principal submatrix $P$ of $A$ containing diagonal
entries $a_{11},a_{22},...,a_{rr},a_{tt},$ we get that
\begin{equation}
\label{bb9}\sum_{i=1}^{r}\lambda_{i}(P)\leq\sum_{i=1}^{r}a_{ii}-\frac{\sum_{s=1}^{r} \left\vert\ a_{ts}\right\vert^{2}}{a_{tt}-\min_{i=1,...,r,t}\left(a_{ii}-\sum_{\substack{s=1\\s\neq i}}^{r+1} \left\vert\ a_{is}\right\vert\right)}. \tag{2.9}
\end{equation}
By the interlacing inequalities \eqref{bb2}, $\sum_{i=1}^{r}\lambda_{i}(A)\leq\sum_{i=1}^{r}\lambda_{i}(P)$. So, \eqref{bb9} gives \eqref{aa9}.   \qed
\subsection{Proof of Theorem 3}
By the Cauchy interlacing principle \eqref{bb2}, the largest eigenvalue of $A$ is greater than or equal to the largest eigenvalue of any $2\times 2$ principal submatrix of $A$. Further, the eigenvalues of  $ \left[ 
		\begin{array}{cc}
			 a_{rr}  & a_{rs}  \\ 
			 \overline{a_{rs}} & a_{ss} %
		\end{array}%
		\right]$are $\frac{1}{2}\left(a_{rr}+a_{ss}\pm\sqrt{(a_{rr}-a_{ss})^{2}+4|a_{rs}|^{2}}\right)$. On using these two facts, we see that
\begin{equation}
\label{bb10}\lambda_{n}(A)\geq a_{nn}+\frac{\sqrt{(a_{nn}-a_{kk})^{2}+4|a_{kn}|^{2}}-(a_{nn}-a_{kk})}{2}  \tag{2.10}
\end{equation}
for all $k=1,2,...,n-1.$
Combining \eqref{bb7} and \eqref{bb10}, we find that 
\begin{equation}
\label{bb11}\sum_{i=1}^{n-1}\lambda_{i}(A)\leq\sum_{i=1}^{n-1}a_{ii}-\frac{\sqrt{(a_{nn}-a_{kk})^{2}+4|a_{kn}|^{2}}-(a_{nn}-a_{kk})}{2}. \tag{2.11}
\end{equation}
Apply \eqref{bb11} to the principal submatrix $Q$ of $A$ containing $
a_{11},a_{22},...,a_{rr},a_{tt}$, we find that for $k=1,2,...,r,$ we have
\begin{equation}
\label{bb12}\sum_{i=1}^{r}\lambda_{i}(Q)\leq\sum_{i=1}^{r}a_{ii}-\frac{\sqrt{(a_{tt}-a_{kk})^{2}+4|a_{tk}|^{2}}-(a_{tt}-a_{kk})}{2}.   \tag{2.12}
\end{equation}
The inequality \eqref{bb12} yields \eqref{aa10}, on using the interlacing inequalities
\eqref{bb2}.    \qed
\par We show by means of an example that \eqref{aa9} and \eqref{aa10} are independent.
\begin{example}
Let 
\begin{equation*}
	   A=\left[ 
		\begin{array}{ccc}
			 2  & 1 & 1  \\ 
			 1  & 2 & 1  \\ 
		 	 1  & 1 & 3 %
		\end{array}%
		\right],\quad
     B=\left[ 
		\begin{array}{ccc}
			 1  & 2 & 3  \\ 
			 2  & 1 & 4  \\ 
		 	 3  & 4 & 1 %
		\end{array}%
		\right].
	\end{equation*}%
Then \eqref{aa9} gives the estimate $\lambda_{1}(A)+\lambda_{2}(A)<\frac{10}{3}$, while \eqref{aa10} gives the weaker estimate $\frac{9-\sqrt{5}}{2}$ for the same quantity. On the other hand from (1.9) we get that  $\lambda_{1}(B)+\lambda_{2}(B)<-\frac{11}{7}$, while from \eqref{aa10} we see that the same quantity is not bigger than $-2$.
\end{example}


\begin{thebibliography}{}

\bibitem{bha97}Bhatia R., \emph{Matrix Analysis}, Springer, New York, (1997).

\bibitem{bha00}Bhatia R., Davis C., \emph{A better bound on the variance}, Amer. Math. Monthly, 107, (2000), 353-357.

\bibitem{bha07}Bhatia R.,  \emph{Positive Definite Matrices}, Princeton University Press, (2007).


\bibitem{hor13}Horn R.A., Johnson C.R., \emph{Matrix Analysis}, Cambridge University Press, (2013).

\bibitem{mar79}Marshal A.W., Olkin I., \emph{Inequalities: Theory of Majorisation and its applications}, Academic Press, (1979).




\end{thebibliography}
\end{document}